\documentclass[12pt,leqno]{article}
\usepackage[pdftex]{graphicx}
\usepackage{amsfonts}

\newcounter{conjecture}
\newcounter{remark}
\newcounter{corollary}

\newenvironment{corollary}{\medskip{\bf Corollary \thecorollary.}
\addtocounter{corollary}{1}\em}{\rm}
\newtheorem{theorem}{Theorem}
\newtheorem{lemma}{Lemma}
\newtheorem{proposition}{Proposition}

\newcommand{\aaa}{\alpha}

\newcommand{\lll}{\label}
\newcommand {\rrr}[1]{(\ref{#1})}

\def \be{\begin{equation}}
\def \ee{\end{equation}}
\def \bt{\begin{theorem}}
\def \et{\end{theorem}}
\def \bc{\begin{corollary}}
\def \ec{\end{corollary}}
\def \bea{\begin{eqnarray}}
\def \eea{\end{eqnarray}}
\def \bas{\begin{eqnarray*}}
\def \eas{\end{eqnarray*}}

\def \noi{\noindent}
\def \aa{\alpha}

\def \ka{\kappa}

\def \si{\sigma}

\def \th{\theta}

\def \pp{\phi}

\def \vski{\vspace{12pt}}

\def \({\left(}
\def \){\right)}

\def \nn{\nonumber}

\def \bc{\begin{center} }
\def \ec{\end{center} }
\def \bs{\begin{slide} }
\def \es{\end{slide} }

\def\square{{\vcenter{\vbox{\hrule height.3pt
         \hbox{\vrule width.3pt height5pt \kern5pt
            \vrule width.3pt}
         \hrule height.3pt}}}}
\def\qed{{\hfill $\square$ \bigskip}}

\newcounter{cccases}
\newcommand{\ccases}[1]{\begingroup \refstepcounter{cccases} {\bf \fontsize{14}{16}\selectfont Case \thecccases }  \label{#1}\endgroup}

\begin{document}

\title{A conjecture of Biggs concerning the resistance of a distance-regular graph.}

\author{
\begin{tabular}{cc}
\textit{Greg Markowsky} & \textit{Jacobus Koolen} \\
gmarkowsky@gmail.com & jacobus\_koolen@yahoo.com \\
Pohang Mathematics Institute & Department of Mathematics\\
POSTECH & POSTECH \\
Pohang, 790-784 & Pohang, 790-784\\
Republic of Korea & Republic of Korea
\end{tabular}}

\maketitle

\begin{abstract}
In \cite{biggs2}, Biggs conjectured that the resistance between any two points on a distance-regular graph of valency greater than 2 is bounded by twice the resistance between adjacent points. We prove this conjecture, give the sharp constant for the inequality, and display the graphs for which the conjecture most nearly fails. Some necessary background material is included, as well as some consequences.
\end{abstract}

\section{Introduction}

The main goal of this paper is to prove the following conjecture of Biggs:

\bt \lll{bigguy}
Let G be a distance-regular graph with degree larger than 2 and diameter D. If $d_j$ is the electric resistance between
any two vertices of distance $j$, then

\be
\max_j d_j = d_D \leq K d_1
\ee

\noi where $K=1 + \frac{94}{101} \approx 1.931$. Equality
holds only in the case of the Biggs-Smith graph.
\et

\noi  We remark that for degree $2$ the theorem is trivially false. This theorem implies several statements concerning random walks on distance-regular graphs, which will
be given at the end of the paper.
General background material on the concept of electric resistance, as well as its connection to random walks, can be found
in the excellent references \cite{doysne} and \cite{biggs}. Biggs' conjecture originally appeared in \cite{biggs2},
which discusses electric resistance on distance-regular graphs only. To understand the proof of the
conjecture, one must understand much of the material in \cite{biggs2}. We have therefore decided to include the
material from
\cite{biggs2} which is key to Theorem \ref{bigguy}. This appears in Section 3, following the relevant graph-theoretic definitions in Section 2. Section 4 gives our proof of the theorem, and Section 5 gives
some consequences, including several in the field of random walks.

\section{Distance-regular graphs}
All the graphs considered in this paper are finite, undirected and
simple (for unexplained terminology and more details, see for example \cite{drgraphs}). Let $G$ be a
connected graph and let $V=V(G)$ be the vertex set of $G$. The distance $d(x,y)$ between
any two vertices $x,y$ of $G$
is the length of a shortest path between $x$ and $y$ in $G$. The diameter of $G$ is the maximal distance
occurring in $G$ and we will denote this by $D = D(G)$.  
For a vertex $x \in V(G)$, define $K_i(x)$ to be the set of
vertices which are at distance $i$ from $x~(0\le i\le
D)$ where $D:=\max\{d(x,y)\mid x,y\in V(G)\}$ is the diameter
of $G$. In addition, define $K_{-1}(x):=\emptyset$ and $K_{D+1}(x)
:= \emptyset$. We write $x\sim_{G} y$ or simply $x\sim y$ if two vertices $x$ and $y$ are adjacent in $G$. A connected graph $G$ with diameter $D$ is called
{\em distance-regular} if there are integers $b_i,c_i$ $(0 \le i
\leq D)$ such that for any two vertices $x,y \in V(G)$ with $d(x,y)=i$, there are precisely $c_i$
neighbors of $y$ in
$K_{i-1}(x)$ and $b_i$ neighbors of $y$ in $K_{i+1}(x)$
(cf. \cite[p.126]{drgraphs}). In particular, distance-regular graph $G$ is regular with valency
$k := b_0$ and we define $a_i:=k-b_i-c_i$ for notational convenience. 
The numbers $a_i$, $b_{i}$ and $c_i~(0\leq i\leq D)$ are called the {\em
intersection numbers} of $G$. Note that $b_D=c_0=a_0=0$, $b_0 = k$ and $c_1=1$.
The intersection numbers of a distance-regular graph $G$ with diameter $D$ and valency $k$ satisfy
(cf. \cite[Proposition 4.1.6]{drgraphs})\\

(i) $k=b_0> b_1\geq \cdots \geq b_{D-1}$;\\
(ii) $1=c_1\leq c_2\leq \cdots \leq c_{D}$;\\
(iii) $b_i\ge c_j$ \mbox{ if }$i+j\le D$.\\

Moreover, if we fix a vertex $x$ of $G$, then $|K_i|$ does not depend on the
choice of $x$ as $c_{i+1} |K_{i+1}| =
b_i |K_i|$ holds for $i =1, 2, \ldots D-1$. In the next section, it will be shown that the resistance between any two vertices
of $G$ can be calculated explicitly using only the intersection array, so that the proof can be conducted using only
the known properties of the array.

\section{Electric resistance on distance-regular graphs}

Henceforth let $G$ be a distance-regular graph with $n$ vertices, degree $k \geq 3$, and diameter $D$. Let $V=V(G)$ and $E=E(G)$ be the vertex and
edge sets, respectively, of $G$. To calculate the resistance between any two vertices we use Ohm's Law, which
states that

\be \label{}
V=IR
\ee

\noi where $V$ represents a difference in voltage(or potential), $I$ represents current, and $R$ represents resistance.
That is, we imagine that our graph is a circuit where each edge is a wire with resistance 1. We attach a battery of
voltage $V$ to
two distinct vertices $u$ and $v$, producing a current through the graph. The resistance between the $u$ and $v$
is then $V$
divided by the current produced. The current flowing through the circuit can be determined by calculating the voltage at
each point on the graph, then summing the currents flowing from $u$, say, to all vertices adjacent to $u$. Calculating the voltage at each point is thereby seen to be an important problem. A function $f$ on $V$ is harmonic at a point $z\in V$ if $f(z)$ is the average of neighboring values of $f$, that is

\be \label{}
\sum_{x\sim z} (f(x)-f(z)) = 0
\ee

The voltage
function on $V$ can be characterized as the unique function which is harmonic on $V-\{u,v\}$ having the prescribed values
on $u$ and $v$. For our purposes, on the distance-regular graph $G$, we will first suppose that $u$ and $v$ are adjacent.
It is easy to see that, for any vertex $z$, $|d(u,z)-d(v,z)| \leq 1$, where $d$ denotes the ordinary graph-theoretic
distance. Thus, any $z$ must be contained in a unique set of one of the following forms:

\bea \label{}
&& K_i^i  = \{x: d(u,x)=i \mbox{ and } d(v,x)=i \}
\\ \nn && K_i^{i+1} = \{x: d(u,x)=i+1 \mbox{ and } d(v,x)=i \}
\\ \nn && K_{i+1}^i = \{x: d(u,x)=i \mbox{ and } d(v,x)=i+1 \}
\eea

\noi Suppose that $(b_0,b_1,\ldots ,b_{D-1};c_1,c_2,\ldots,c_D)$ is the intersection array of $G$. For $0 \leq i \leq D-1$
define the numbers $\phi_i$ recursively by

\bea \label{smile}
&& \pp_0=n-1
\\ \nn && \pp_i = \frac{c_i\pp_{i-1}-k}{b_i}
\eea

\noi We then have the following fundamental proposition.

\begin{proposition} \label{vp}
The function $f$ defined on $V$ by

\bea \label{}
&& f(u) = -f(v) = \pp_0
\\ \nn && f(z)= 0 \mbox{ for } x \in K_{i}^i
\\ \nn && f(z)= \pp_i \mbox{ for } x \in K_{i+1}^i
\\ \nn && f(z)= -\pp_i \mbox{ for } x \in K^{i+1}_i
\eea

is harmonic on $V-\{u,v\}$.
\end{proposition}

\noi In the following intersection diagram, the value of $f$ on each set is given directly
outside the set.

\includegraphics[width=130mm,height=110mm]{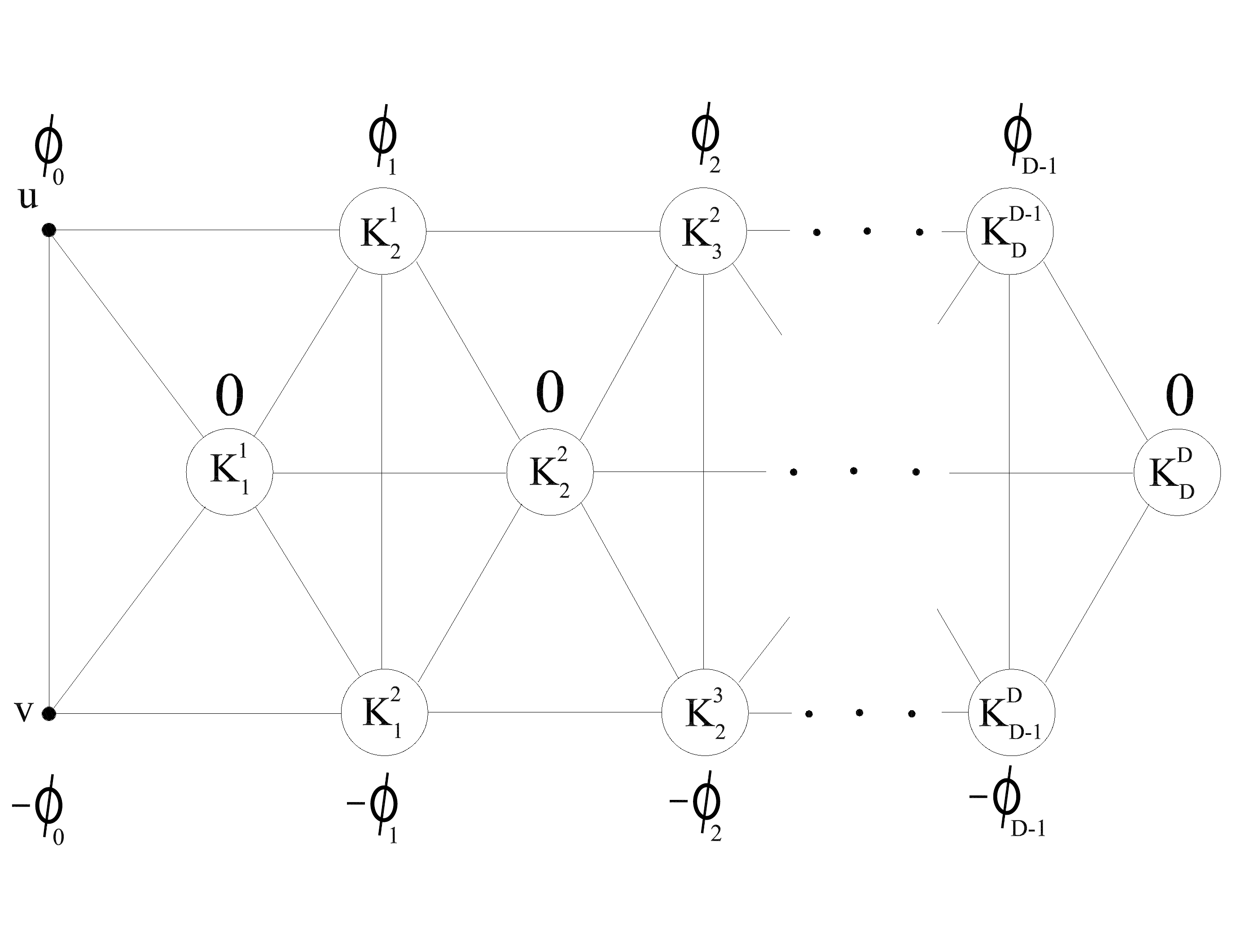}
\begin{center} Figure 1 \end{center}

To prove Proposition \ref{vp} we need the following lemma, which may be of interest in its own right.

\begin{lemma} \label{hh}
Let $z \in G$, and let $K_i = \{x:d(z,x)=i\}$ as in Section 2. Let $e_i$ be the number of edges of $G$ with one endpoint in $K_{i}$
and the other in $K_{i+1}$. Then
\be \lll{sam}
\pp_i = \frac{k \sum_{j > i}|K_j|}{e_i}
\ee
\end{lemma}

\noi {\bf Proof:} Since $\pp_0=n-1=\sum_{j > 0}|K_j|$ and $e_0=k$, it is clear that \rrr{sam} holds for $i=0$. We need therefore
only verify that the numbers $\psi_i = \frac{k \sum_{j > i}|K_j|}{e_i}$
satisfy the recursive relation given in \rrr{smile}. This is immediate from the facts that $e_i = b_i|K_i|$ and
$e_{i-1} = c_i|K_i|$, for we see that

\bea \label{}
&& \frac{c_i \psi_{i-1} - k}{b_i} = \frac{c_i(\frac{k |K_i| + k \sum_{j > i}|K_j|}{e_{i-1}}) - k}{b_i}
\\ \nn && \hspace{2cm} = \frac{c_i k \sum_{j > i}|K_j|}{b_i e_{i-1}}
\\ \nn && \hspace{2cm} = \frac{k \sum_{j > i}|K_j|}{b_i K_i}
\\ \nn && \hspace{2cm} = \frac{k \sum_{j > i}|K_j|}{e_i}
\\ \nn && \hspace{2cm} = \psi_i
\eea \qed

\noi {\bf Proof of Proposition \ref{vp}:} Suppose first that $z \in K_i^i$ for some $i$. The points adjacent to $z$ must
lie within
$K_i^i \bigcup K_{i-1}^{i-1} \bigcup K_{i+1}^{i+1} \bigcup K_{i}^{i+1} \bigcup K_{i}^{i-1} \bigcup K_{i-1}^{i}
\bigcup K_{i+1}^{i}$. Since $b_i$ is equal to the number of adjacent points in $K_{i+1}^{i+1} \bigcup K_{i}^{i+1}$, and
also in the set $K_{i+1}^{i+1} \bigcup K_{i+1}^{i}$, we see that

\be
|\{x:z \sim x \mbox{ and } x \in K_{i+1}^{i}\}| = |\{x:z \sim x \mbox{ and } x \in K_{i}^{i+1}\}|
\ee

\noi A similar argument shows

\be
|\{x:z \sim x \mbox{ and } x \in K_{i-1}^{i}\}| = |\{x:z \sim x \mbox{ and } x \in K_{i}^{i-1}\}|
\ee

\noi It follows from this that

\be \label{}
\sum_{x\sim z}f(x) = 0 = f(z)
\ee

\noi and $f$ is harmonic at $z$. Now suppose that $z \in K^{i}_{i+1}$ with $1 \leq i \leq D-2$. Here
the points adjacent to $z$ must
lie within
$K^i_{i+1} \bigcup K^{i-1}_{i} \bigcup K^{i+1}_{i+2} \bigcup K_{i}^{i} \bigcup K_{i+1}^{i+1} \bigcup K^{i+1}_{i}$.
The number of edges from $z$ to points in $K^{i-1}_{i}$ is $c_i$ and to points in $K^{i+1}_{i+2}$ is $b_{i+1}$. Let
the number of edges from $z$ to points in $K^{i+1}_{i}$ be $\aaa$. Then the number of edges from $z$
to other points in $K^i_{i+1}$ is given by $k+\aa -c_{i+1}-b_i$. We therefore have

\bea \label{}
&& \sum_{x\sim z} f(x) = b_{i+1} \pp_{i+1} + c_i \pp_{i-1} + (k+\aa -c_{i+1}-b_i)\pp_i + \aa (-\pp_i)
\\ \nn && \hspace{1.6cm} = k \pp_i = k f(z)
\eea

\noi where we have used the following equations equivalent to the recursive relation in \rrr{smile}.

\bea \label{pt}
&& c_i \pp_{i-1} = b_i\pp_i +k
\\ \nn && b_{i+1} \pp_{i+1} = c_{i+1}\pp_i - k
\eea

\noi We see that $f$ is harmonic at $z$. The same argument works for $z \in K^{D-1}_{D}$, except that there is some
difficulty in using the last equation in \rrr{pt}, as $b_D=0$, and $\pp_i$ was only defined for $i \leq D-1$.
Happily, Lemma \ref{hh} solves our dilemma, for as an immediate consequence we obtain $\pp_{D-1} = \frac{k}{c_D}$.
Thus, defining $\pp_D = 0$ is consistent with \rrr{pt}, and $f$ is harmonic on $K^{D-1}_{D}$. By symmetry, $f$ is harmonic
at all points lying in sets of the form $K^{i}_{i+1}$, and the proof is complete.
\qed

\begin{corollary}
$\phi_i > \phi_{i+1}$ for $0 \leq i \leq D-2$
\end{corollary}

\vski

{\bf Proof:} Suppose $\phi_i \leq \phi_{i+1}$ for some $i$. Due to the monotonicity of the sequences $b_i,c_i$, we would have

\be \label{}
\phi_{i+2} = \frac{c_{i+2}\phi_{i+1}-k}{b_{i+2}} \geq \frac{c_{i+1}\phi_{i}-k}{b_{i+1}} = \phi_{i+1}
\ee

Continuing in this way we would have $\phi_{D-1} \geq \phi_{D-2}$. On the other hand, by harmonicity $\phi_{D-1}$ is the weighted average of the values $\phi_{D-2},0,$ and $-\phi_{D-1}$, so that $\phi_{D-1} < \phi_{D-2}$. This is a contradiction. \qed

\noi It may interest the reader to note that the subtracted constant $k$ in the numerator of the recursive relation of
\rrr{smile} can be replaced by any constant without affecting harmonicity outside of the sets $K^{D-1}_{D}$ and
$K^{D}_{D-1}$.
However, $k$ is the only constant which gives $\pp_D = 0$, and therefore is the constant dictated by the requirement
that $f$
be harmonic and attain the boundary values of $(n-1)$ and $-(n-1)$ at $u$ and $v$. The resistance between $u$ and $v$ can
now easily be computed as the voltage difference between the points, $2\pp_0=2(n-1)$, divided by the current $I$
flowing through
the circuit. This current is the sum of the voltage differences between $u$ and vertices adjacent to $u$, and is
readily computable as $I=nk$. We see that the resistance between $u$ and $v$ is

\be \lll{calc}
R_{uv}=\frac{2(n-1)}{nk}=\frac{n-1}{m}
\ee

\noi where $m=nk/2$ is the number of edges in $G$. This result is in fact an immediate consequence of Foster's Network Theorem(see \cite{biggs} or \cite{foster}), and was derived, among other things, by other methods in \cite{dkb}. In the remainder of this section, however, it will be more conceptually convenient to keep
$I$ and the $\pp$'s
in the formulas rather than their explicit values, as this reminds us that they represent the current and voltages,
respectively. Calculating the resistances between nonadjacent vertices might now seem to be a formidable task, but in
fact there is virtually no more to be done. We have the following proposition.

\begin{proposition} \label{}
The resistance between two vertices of distance $j$ in a graph is given by
\be \label{}
\frac{2\sum_{0\leq i<j}\pp_i}{I}
\ee
\end{proposition}

\noi {\bf Proof:} Suppose $d(u,v)=j$. We can choose points $x_0=u, x_1, \ldots , x_j = v$ such that $x_i \sim x_{i+1}$. For
any pair of adjacent points $y,z$ we let $f_{yz}$ be the unique function on $V$ given in Proposition \ref{vp} which is
harmonic on $V-\{y,z\}$ and which satisfies $f(w)=-f(z)=\pp_0$. The key claim is that for any three points
$w,y,z$ with $y \sim w \sim z$ the function $f_{yw}+f_{wz}$ is harmonic on $V - \{y,z\}$. This is clear for all points
in $V - \{y,z\}$ except $w$. To show harmonicity at $w$, note that a current of $I$ flows into $w$ due to $f_{yw}$, whereas a current of $I$ flows out of $w$ due to $f_{wz}$. The net current flow into $w$ is therefore $0$, which is
equivalent to harmonicity(see \cite{doysne}). Thus, the voltage function
$g=\sum_{0\leq i \leq j-1}f_{x_ix_{i+1}}$, which is harmonic on $V-\{u,v\}$, gives rise to a current of $I$ flowing from
$u$ to $v$. We must therefore calculate the values of the function $g$ at the points $u$ and $v$. It is straightforward
to verify that $f_{x_ix_{i+1}}(u) = \pp_i$(since $u$ lies in the set $K_{i+1}^i$ formed with respect
to the pair $x_i,x_{i+1}$), and likewise $f_{x_ix_{i+1}}(v) = -\pp_{D-(i+1)}$. Thus, $g(u)=\sum_{0\leq i<j}\pp_i$ and
$g(v)=-\sum_{0\leq i<j}\pp_i$. The result follows. \qed

\section{Proof of Theorem}

In fact, we will prove a statement stronger than Theorem \ref{bigguy}. Let $\mathcal{E}$ be the set of the
following four graphs, with corresponding properties listed:

\vski

\noi \begin{tabular}{ l c c r }
Name\footnotemark & Vertices & Intersection array & $\frac{\phi_1+ \ldots + \phi_{D-1}}{\phi_0}$ \\
\hline
Biggs-Smith Graph & 102 & (3,2,2,2,1,1,1;1,1,1,1,1,1,3) & 0.930693\\
Foster Graph & 90 & (3,2,2,2,2,1,1,1;1,1,1,1,2,2,2,3) & 0.896067\\
Flag graph of $GH$(2,2) & 189 & (4,2,2,2,2,2;1,1,1,1,1,2) & 0.882979\\
Tutte's 12-Cage & 126 & (3,2,2,2,2,2;1,1,1,1,1,3) & 0.872
\end{tabular}

\footnotetext{The referee has pointed out that Tutte's 12-Cage may be more accurately referred to as Benson's graph, and indeed the literature is mixed on this point. The referee further remarked that the Flag graph of $GH$(2,2) can also be realized as the line graph of Tutte's 12-Cage, or Benson's graph. In this table, we are employing the names given in \cite{drgraphs}.}

\bt \lll{db}
Other than graphs in $\mathcal{E}$, for any distance regular graph with degree at least 3 we have
\be
\phi_1+ \ldots + \phi_{D-1} < .87 \phi_0
\ee
\et

\noi This clearly implies Theorem \ref{bigguy} and shows that the graphs in $\mathcal{E}$ are the extremal cases.

\vski

\noi {\bf Proof of Theorem \ref{db}:} The proof proceeds by considering a number of separate cases, and leans heavily on the standard reference \cite{drgraphs}. Without access to this book, the proof will likely be incomprehensible to the reader. In the
estimates used in the proof, the $-k$ in the numerator of the recurrence relation is largely ignored, but the reader
should be warned that this term is by no means unnecessary.
That is because it is crucial that the $\phi_i$'s form a monotone decreasing sequence, and without the $-k$ this would
not be the case. Nevertheless, we will from this point forth mainly use the facts $\phi_i < \frac{c_i \phi_{i-1}}{b_i}$
and $\phi_i < \phi_{i-1}$. We are required to show

\be
\frac{\phi_1+ \ldots + \phi_{D-1}}{\phi_0} \leq .87
\ee

\noi for all graphs not in $\mathcal{E}$.

\vski

\ccases{tit}: $D = 2$.

We need only show $\phi_1<.87\phi_0$. This is clear if $b_1>1$, since $c_1=1$ and $\phi_i < \frac{c_i \phi_{i-1}}{b_i}$. The case $b_1=1$ is known to occur only in the case of the Cocktail party graphs, and it is simple to verify the relation in this case. \qed

\ccases{deg3}: $k=3$.

It is known(see \cite{drgraphs}, Theorem 7.5.1) that the only distance-regular graphs of degree 3 with
diameter greater than 2 are given by the intersection arrays below, and which give rise to the resistances given:

\begin{tabular}{ l c c r }
Name & Vertices & Intersection array & $\frac{\phi_1+ \ldots \phi_{D-1}}{\phi_0}$ \\
\hline
Cube & 8 & (3,2,1;1,2,3)  &   0.428571 \\
Heawood graph & 14 & (3,2,2;1,1,3)  &  0.461538 \\
Pappus graph & 18 & (3,2,2,1;1,1,2,3)  &      0.588235 \\
Coxeter graph & 28 & (3,2,2,1;1,1,1,2)   &     0.666667 \\
Tutte's 8-cage & 30 & (3,2,2,2;1,1,1,3)    &    0.655172 \\
Dodecahedron & 20 & (3,2,1,1,1;1,1,1,2,3) &   0.842105 \\
Desargues graph & 20 & (3,2,2,1,1;1,1,2,2,3)  &  0.710526 \\
Tutte's 12-cage & 126 & (3,2,2,2,2,2;1,1,1,1,1,3)  &     0.872 \\
Biggs-Smith graph & 102 & (3,2,2,2,1,1,1;1,1,1,1,1,1,3) &  0.930693 \\
Foster graph & 90 & (3,2,2,2,2,1,1,1;1,1,1,1,2,2,2,3)  &      0.896067 \\
\end{tabular}

\qed

\ccases{deg4}: $k=4$.

It is known(see \cite{broukool}) that the only distance-regular graphs of degree 4 with diameter greater than 2
are given by the intersection arrays below, and which give rise to the resistances given:

\noi \begin{tabular}{ l c c r }
Name & Vertices & Intersection array & $\frac{\phi_1+ \ldots \phi_{D-1}}{\phi_0}$ \\
\hline
$K_{5,5}$ minus a matching& 10 & (4,3,1;1,3,4) &   0.296296\\
Nonincidence graph of $PG(2,2)$ &  14 & (4,3,2;1,2,4)  &  0.307692\\
Line graph of Petersen graph &  15 & (4,2,1;1,1,4)  &  0.428571\\
4-cube &  16 & (4,3,2,1;1,2,3,4) &       0.422222\\
Flag graph of $PG(2,2)$ & 21 & (4,2,2;1,1,2) &   0.5\\
Incidence graph of $PG(2,3)$ & 26 & (4,3,3;1,1,4)  &  0.32\\
Incidence graph of $AG(2,4)$-p.c. & 32 & (4,3,3,1;1,1,3,4)  &      0.376344\\
Odd graph $O_4$ & 35 & (4,3,3;1,1,2)  &  0.352941\\
Flag graph of $GQ(2,2)$ & 45 & (4,2,2,2;1,1,1,2)  &      0.681818\\
Doubled odd graph & 70& (4,3,3,2,2,1,1;1,1,2,2,3,3,4)&    0.521739\\
Incidence graph of $GQ(3,3)$ & 80 & (4,3,3,3;1,1,1,4) &       0.417722\\
Flag graph of $GH(2,2)$ & 189 & (4,2,2,2,2,2;1,1,1,1,1,2) &      0.882979\\
Incidence graph of $GH(3,3)$ & 728 & (4,3,3,3,3,3;1,1,1,1,1,4) &      0.485557\\
\end{tabular} \qed

\ccases{mzy}: $D \leq 5, b_1 \geq 5$.

This case was done initially by Biggs in \cite{biggs2}, without the restriction on $b_1$ but with the constant $1$ in place of $.87$. Nevertheless, when we restrict $b_1$ as above this is trivial, because $\frac{\phi_1}{\phi_0}<\frac{1}{b_1}$ and $\phi_i \leq \phi_1$ for all $i>0$. Therefore,

\be
\frac{\phi_1+ \ldots + \phi_{D-1}}{\phi_0} \leq \frac{(D-1)\phi_1}{\phi_0} \leq \frac{4}{b_1} \leq .8
\ee \qed

Henceforth, in all cases for which $b_1 \geq 5$ we can assume $D \geq 6$. In what follows, let $j$ denote the smallest value such
that $c_j \geq b_j$. If $c_j > b_j$, then, since $c_{D-j} \leq b_j$ and the $c_i$'s are nondecreasing, we see that
$D-j<j$, hence $D \leq 2j-1$. If $c_j=b_j$, then it follows from Corollary 5.9.6 of \cite{drgraphs} that $c_{2j}>b_{2j}$.
For this to occur, either $c_{2j}>b_j$ or $c_j > b_{2j}$. By the same argument as before, we obtain $D \leq 3j-1$.
This will be of fundamental importance in our proof. To begin with, we see that when $D \geq 6$ we must have $j \geq 3$.

\vski

\ccases{line}: $G$ is a line graph.

The distance-regular line graphs have been classified, and appear in Theorem 4.2.16 of \cite{drgraphs}. All such
graphs with $k \geq 3$ have $D \leq 2$ and are therefore covered by Case \ref{tit}, with two exceptions. First of all, $G$ may be a generalized $2D$-gon of order $(1,s)$. The intersection array of $G$ is then of the form $(2(a_1+1),a_1+1, \ldots, a_1+1;1,1,\ldots,1,2)$, with $a_1>1$. The other possibility is that $G$ could be the line graph of a Moore graph, and in this case the intersection array of $G$ is of the form $(2\ka -2, \ka - 1, \ka - 2; 1,1,4)$, for some $\ka \geq 3$. In both of these cases it is straightforward to verify that the conclusion of the theorem holds. \qed

\vski

\ccases{cm}: $b_1 \geq 5, j=3, c_2=1$.

Since $j=3$, $b_2 \geq 2$ and $D \leq 8$. We have

\be \label{}
\frac{\phi_1+ \ldots + \phi_{D-1}}{\phi_0} \leq \frac{\phi_1+ 6 \phi_{2}}{\phi_0} \leq \frac{1}{b_1}+ \frac{6}{2b_1} = \frac{4}{b_1} \leq .8
\ee \qed

\ccases{jr}: $b_1 \geq 5, j=3, c_2>1$.

By Theorem 5.4.1 in \cite{drgraphs}, $c_2 \leq \frac{2}{3} c_3$. If $c_3 > b_3$ then $D \leq 2j-1 = 5$, which was covered in Case \ref{mzy}. If $c_3=b_3 \leq b_2$, then if we assume $\frac{c_2}{b_2} \leq \frac{1}{2}$ we have

\be \label{}
\frac{\phi_1+ \ldots + \phi_{D-1}}{\phi_0} \leq \frac{\phi_1+ 6\phi_{2}}{\phi_0} \leq \frac{1}{b_1}+ \frac{3}{b_1} = \frac{4}{b_1} \leq .8
\ee

On the other hand, if it is not the case that $\frac{c_2}{b_2} \leq \frac{1}{2}$, then the proof of Theorem 5.4.1 of \cite{drgraphs} implies that $G$ contains a quadrangle. By Corollary 5.2.2 in \cite{drgraphs}, $D\leq\frac{2k}{k+1-b_1}$. It is straightforward to verify that the fact that $k \geq b_1+1$ implies that

\be \label{}
\frac{2k}{k+1-b_1} \leq b_1+1
\ee

We therefore see that the fact that $G$ contains a quadrangle implies $D \leq b_1+1$. Furthermore, we still have $\frac{c_2}{b_2} \leq \frac{2}{3}$ by Theorem
5.4.1 of \cite{drgraphs}. We therefore have

\be \label{}
\frac{\phi_1+ \ldots + \phi_{D-1}}{\phi_0} \leq \frac{\phi_1+ (b_1-1)\phi_{2}}{\phi_0} \leq \frac{1}{b_1}+ \frac{2(b_1-1)}{3b_1} = \frac{2b_1+1}{3b_1} \leq .7
\ee \qed

\ccases{fade}: $b_1 \geq 5, j \geq 4, c_2 = 1$.

If $j \geq 4$ and $b_2=2$ then we must have $b_3=2, c_3=1$, so that $\frac{b_2 b_3}{c_2 c_3} = 4$. On the other hand, if this does not occur than $\frac{b_2}{c_2} \geq 3$. We will consider these cases separately.

\vski

{\bf Subcase 1:} $\frac{b_2}{c_2} \geq 3$.

For $i < j$ we have $b_1 \geq b_i > c_i$, and for any $i$ with $c_i>1$ we must have $b_i<b_1$, by Proposition 5.4.4 in \cite{drgraphs}. Thus, $\frac{c_i}{b_i} \leq \frac{b_1-2}{b_1-1}$. Define $\aa = \frac{b_1-2}{b_1-1}$. We have

\be \label{}
\frac{\phi_1+ \ldots + \phi_{D-1}}{\phi_0} \leq \frac{1}{b_1}+ \frac{1}{3b_1}+ \frac{\aa}{3b_1}+ \ldots + \frac{\aa^{j-3}}{3b_1} + \frac{(2j-1)\aa^{j-3}}{3b_1}
\ee

Replace the second through $(j-1)$th term by a geometric series to obtain

\bea \label{kyl}
\frac{\phi_1+ \ldots + \phi_{D-1}}{\phi_0} < \frac{1}{b_1}+ \frac{1}{3b_1}\Big(\frac{1}{1-\frac{{b_1-2}}{b_1-1}}\Big) + \frac{(2j-1)\aa^{j-3}}{3b_1}
\\ \nn < \frac{1}{b_1} + \frac{b_1-1}{3b_1} + \frac{2(j-1/2)\aa^{j-1/2}}{3b_1\aa^{5/2}}
\eea
Simple calculus shows that the maximum of the function $u\aa^u$ is $\frac{-1}{e \ln{\aa}}$. We therefore obtain

\be \label{hhh}
\frac{\phi_1+ \ldots + \phi_{D-1}}{\phi_0} < \frac{b_1+2}{3b_1} + \frac{-2}{3b_1(\frac{b_1-2}{b_1-1})^{5/2}e\ln(\frac{b_1-2}{b_1-1})}
\ee

It is straightforward to verify that the function $(b-2)\ln(\frac{b-2}{b-1})$ is increasing in $b$, so that the right hand side of \rrr{hhh} achieves its maximum on the allowed range when $b_1=5$. Plugging in $b_1=5$ gives approximately $.851$ as a bound for \rrr{hhh}. \qed

{\bf Subcase 2:} $\frac{b_2 b_3}{c_2 c_3} \geq 4$.

This follows much as in the previous case, except that we may simplify by using the slightly weaker bound $\frac{c_i}{b_i} \leq \frac{b_1-1}{b_1}$ for $i<j$. Let $\aa = \frac{b_1-1}{b_1}$. Since $b_2 \geq b_3$ and $c_2 \leq c_3$ we must have $\frac{b_2}{c_2} \geq 2$. We then have

\be \label{}
\frac{\phi_1+ \ldots + \phi_{D-1}}{\phi_0} \leq \frac{1}{b_1}\!+ \! \frac{1}{2b_1}\!+ \! \frac{1}{4b_1}\!+ \! \frac{\aa}{4b_1}\!+ \! \ldots \!+ \! \frac{\aa^{j-3}}{4b_1}\! + \!\frac{(2j-1)\aa^{j-3}}{4b_1}
\ee

Following the steps in \rrr{kyl} above, we obtain

\be \label{min}
\frac{\phi_1+ \ldots + \phi_{D-1}}{\phi_0} < \frac{3}{2b_1} + \frac{1}{4} + \frac{-1}{2b_1(\frac{b_1-1}{b_1})^{5/2}e\ln(\frac{b_1-1}{b_1})}
\ee

Again this is decreasing in $b_1$, and plugging in $b_1=5$ gives a bound for \rrr{min} of about $.84$. \qed

\ccases{the}: $b_1 \geq 3, j \geq 4, c_2 > 1, G$ contains a quadrangle.

As in the argument given in Case \ref{jr}, we see that $G$ containing a quadrangle implies $D \leq b_1+1$. Furthermore,
Theorem 5.4.1 of \cite{drgraphs} implies that $c_3 \geq (3/2)c_2$. Since $j \geq 4$ and thus $b_2 \geq b_3 > c_3$
we must have $\frac{c_2}{b_2} \leq \frac{2}{3}$. This gives

\be \label{}
\frac{\phi_1+ \ldots + \phi_{D-1}}{\phi_0} \leq \frac{1}{b_1}+ (b_1-1)\frac{2}{3b_1} = \frac{2b_1+1}{3b_1}
\ee

When $b_1 \geq 3$ this is bounded by $.8$. \qed

\ccases{guy}: $b_1 \geq 3, j \geq 4, c_2 \geq 1, G$ does not contain a quadrangle.

In this case $G$ is a Terwilliger graph. By Corollary 1.16.6 of \cite{drgraphs}, if $k<50(c_2-1)$ then $D \leq 4$ and $b_1 \geq 5$, which was covered in Case \ref{mzy}. Thus, we can assume $k \geq 50(c_2-1)$, which implies $b_1 \geq 10c_2$. If $b_2 \geq 3c_2$ then we can follow the proof of Subcase 1 of Case \ref{fade} to obtain our result, so we may assume $b_2 \leq 3c_2$, which implies $b_2 < \frac{b}{2}$. It follows from this that for $i<j$ $\frac{c_2}{b_2} \leq \frac{(b_1/2)-1}{b_1/2}=\frac{b_1-2}{b_1}$. We set $\aa = \frac{b_1-2}{b_1}$. By the proof of Theorem 5.4.1 in \cite{drgraphs} we have $c_3 \geq 2 c_2$. Since $b_2 \geq b_3 > c_3 \geq 2c_2$ we have $\frac{b_2}{c_2} \geq 2$. We compute

\be \label{}
\frac{\phi_1+ \ldots + \phi_{D-1}}{\phi_0} \leq \frac{1}{b_1}+ \frac{1}{2b_1}+ \frac{\aa}{2b_1}+ \ldots + \frac{\aa^{j-3}}{2b_1} + \frac{(2j-1)\aa^{j-3}}{2b_1}
\ee

Replace the second through $(j-1)$th term by a geometric series to obtain

\bea \label{kyl}
\frac{\phi_1+ \ldots \phi_{D-1}}{\phi_0} < \frac{1}{b_1}+ \frac{1}{2b_1}\Big(\frac{1}{1-\frac{{b_1-2}}{b_1}}\Big) + \frac{(2j-1)\aa^{j-3}}{2b_1}
\\ \nn < \frac{1}{b_1} + \frac{1}{4} + \frac{(j-1/2)\aa^{j-1/2}}{b_1\aa^{5/2}}
\eea
The maximum of the function $u\aa^u$ is $\frac{-1}{e \ln{\aa}}$. We therefore obtain

\be \label{hh7}
\frac{\phi_1+ \ldots \phi_{D-1}}{\phi_0} < \frac{1}{b_1} + \frac{1}{4} + \frac{-1}{b_1(\frac{b_1-2}{b_1})^{5/2}e\ln(\frac{b_1-2}{b_1})}
\ee

As before, the function $(b-2)\ln(\frac{b-2}{b})$ is increasing in $b$, so the right hand side of \rrr{hh7} is decreasing in $b_1$. Plugging in $b_1=10$(recall that $b_1 \geq 10c_2 \geq 10$) gives approximately $.64$ as a bound. \qed

\ccases: $b_1=3$ or $4, k \geq 5, c_2 = 1$.

This will be broken down into cases by degree $k$. Proposition 1.2.1 in \cite{drgraphs} implies that $(a_1+1) | k$, so since $b_1=k-a_1 -1$ and $b_1>0$ we see that $b_1 \geq k/2$. This implies $k \leq 8$.

\vski

{\bf Subcase $k=8$}: $b_1=3$ is ruled out because $(a_1+1) | k$. Suppose $b_1=4$. By Proposition 4.3.4 of \cite{drgraphs}, $G$ is a line graph, and is therefore covered by Case \ref{line}.

\vski

{\bf Subcase $k=7$}: Since $(a_1+1) | k$, we must have $a_1=0$ and thus $b_1=6$, which is a contradiction.

\vski

{\bf Subcase $k=6$}: Since $(a_1+1) | k$, we have $a_1 \in \{0,1,2\}$. If $a_1 = 0$, then $b_1=5$, a contradiction. If $a_1=1$, then as was shown in \cite{hns} $G$ is one of the following graphs.

\begin{tabular}{ l c c r }
Name & Vertices & Intersection array & $\frac{\phi_1+ \ldots \phi_{D-1}}{\phi_0}$ \\
\hline
Colinearity graph of $GQ(2,2)$ & 15 & (6,4;1,3) & 0.142857 \\
Colinearity graph of $GH(2,2)$ & 27 & (6,4,2;1,2,3) & 0.269231 \\
Hamming graph $H(3,3)$ & 63 & (6,4,4;1,1,3) & 0.258065 \\
Halved Foster graph & 45 & (6,4,2,1;1,1,4,6) & 0.278409 \\
\end{tabular}

\vski

If $a_1=2$, then by Proposition 4.3.4 of \cite{drgraphs}, $G$ is a line graph, and
is therefore covered by Case \ref{line}.

\vski

{\bf Subcase $k=5$}: Since $(a_1+1) | k$, we must have $a_1=0$ and $b_1=4$. Suppose first that $b_2=3$ or $4$. Note that, for $i<j$, $\frac{c_i}{b_i} \leq \frac{2}{3}$, since $c_i+b_i \leq 5$. Using the same technique as in many of the previous cases we have

\bea \label{hh8}
&& \frac{\phi_1+ \ldots \phi_{D-1}}{\phi_0} < \frac{1}{4} + \frac{1}{12} + \frac{1}{12}\big(\frac{2}{3}+ \ldots + (\frac{2}{3})^{j-2}\big) + \frac{1}{12}(\frac{2}{3})^{j-2}(2j-1)
\\ \nn && < \frac{1}{4} + \frac{1}{12} + \frac{3}{12} + \frac{1}{16}(\frac{2}{3})^{j-2}(2j-1)
\eea

It is straightforward to verify that the last expression in \rrr{hh8} is decreasing in $j$ for $j \geq 3$. Plugging in $j=3$ gives a bound of $31/36<.87$. It remains only to consider $b_2 \leq 2$. Suppose $b_2=2$. If $c_3=1$, it would follow from Corollary 4.3.12(ii) that $3$ divides $20$. Thus, we can assume $c_3 \geq 2$, and therefore $j=3$ and $D \leq 8$. We will first show that $n \leq 140$. Fix a point $u$ in $G$ and let $k_i = |\{v:d(u,v)=i\}|$. The numbers $k_i$ are easily computable through the intersection arrays by $k_i=\frac{\prod_{l=0}^{i-1}b_i}{\prod_{l=1}^{i}c_i}$. The $k_i$'s are nonincreasing for $i\geq j$, so since $k_3=20$, if $D\leq 7$ we have $n \leq 1+5+6(20) < 140$. Suppose $D=8$. Then $c_6 \leq b_2 =2$, so $c_6=2$ and this implies $b_6=1$. In this case, $k_7 =10$, and thus $k_8 \leq 10$ as well. We get $n \leq 1+5+5(20) + 2(10) < 140$ again. Since $k=5$, we get $k>(n-1)/28$. Let $\th=|\{i:b_i=c_i=2\}|$. If $\th=3$, the maximal allowed value, we have the following calculations:

\bea \label{} \nn
&& \phi_0 = n-1 , \phi_1 < \frac{n-1}{4}, \phi_2 < \frac{n-1}{8} ,
\\ \nn && \phi_3 < \frac{2((n-1)/8) - (n-1)/28}{2} = \frac{6(n-1)}{56},
\\ \nn && \phi_4 < \frac{2(6(n-1)/56) - (n-1)/28}{2} = \frac{5(n-1)}{56},
\\ \nn && \phi_5 < \frac{2(5(n-1)/56) - (n-1)/28}{2} = \frac{4(n-1)}{56}
\eea

Since $\phi_6,\phi_7 < \phi_5$ we get

\be
\frac{\phi_1+ \ldots \phi_{D-1}}{\phi_0} < \frac{1}{4} + \frac{1}{8} + \frac{6}{56} + \frac{5}{56} + 3\Big(\frac{4}{56}\Big) = \frac{44}{56} < .87
\ee

Similar but easier calculations handle the cases $\th=2,1,0$. The case $b_2=1$ can also be handled in a similar way. Note that in this case $j=2$, so $D \leq 5$. If $D \leq 4$, then

\be \label{}
\frac{\phi_1+ \ldots \phi_{D-1}}{\phi_0} \leq \frac{3 \pp_1}{\pp_0} < \frac{3}{4}
\ee

If $D=5$, then $k_1=5$, $k_2=20$, and $k_i \leq 20$ for $i\geq j$(since the $k_i$'s are nonincreasing for $i\geq j$). It follows that $n\leq 86$, and therefore $k>\frac{n-1}{20}$. Furthermore, $c_3 \leq b_2 =1$, so $c_3=1$. Thus,

\bea \label{} \nn
&& \phi_0 = n-1 , \phi_1 < \frac{(n-1)-k}{4} \leq \frac{19(n-1)}{80},
\\ \nn && \phi_2 < \frac{c_2\phi_1-k}{1}<\frac{15(n-1)}{80},
\\ \nn && \phi_3,\pp_4 < \frac{c_2\phi_2-k}{1}<\frac{11(n-1)}{80}
\eea

And so

\be \label{}
\frac{\phi_1+ \ldots \phi_{4}}{\phi_0} < \frac{19}{80}+\frac{15}{80}+2\left(\frac{11}{80}\right) =56/80 =.7
\ee

\qed

\section{Consequences}

As indicated in \cite{biggs2}, there are some immediate consequences for random walks. Let $u$ be a vertex of $G$, and
and suppose we start a random walk at $u$. For any other point $v$, we let the expected number of steps needed to hit
$v$ be denoted $H_{uv}$. This is referred to as the {\it hitting time}. The {\it commute time} $C_{uv}$ is the expected
number of
steps necessary for the random walk to travel from $u$ to $v$ and back to $v$, and in the case of distance regular graphs
is equal to $2H_uv$. By Theorem 1 in \cite{comcov},
the expected commute time of a random walk between two points $u$ and $v$ is equal to $2mR_{uv}$. Thus, from Theorem
\ref{bigguy} in this paper, and the calculation of resistance given in Section 2, in a distance-regular graph with valency greater than 2 we have

\begin{proposition} \label{}
\be
H_{uv} \leq 2m\Big( \frac{n-1}{m} \Big) = 2(n-1)
\ee
\be
C_{uv} \leq 4m\Big( \frac{n-1}{m} \Big) = 4(n-1)
\ee
\end{proposition}

\noi The {\it cover time} $Co(G)$ is the expected number of steps that our random walk requires before it has visited every
site on $G$. Applying Theorem 3 in \cite{comcov}, we have

\begin{proposition} \label{tree}
For $n$ large,

\be \label{}
Co(G) \leq (4+o(1))(n-1)\ln{n}
\ee
\end{proposition}

\noi In fact, in \cite{fue} it was shown that for all graphs, distance-regular or otherwise, we have

\be \label{}
Co(G) \geq (1+o(1))n\ln{n}
\ee

so that the bound in Proposition \ref{tree} is the best possible, up to the multiplicative constant. Let $\si$ be the smallest nonzero eigenvalue of the Laplacian matrix. Note that $k-\si$ is the second largest eigenvalue of the adjacency matrix. Let $R_{max}$ denote the largest resistance
between points in $G$, which we have seen necessarily occurs when the points are at distance $D$. Combining Theorem
\ref{bigguy} in this
paper with Theorem 7 in \cite{comcov}, we have

\begin{proposition} \label{}
\be \label{}
\si \geq \frac{1}{n R_{max}} \geq \frac{m}{2n(n-1)} = \frac{k}{4(n-1)}
\ee
\end{proposition}

There have been discussions between the two authors as to whether Theorem \ref{bigguy} really gives new information on the structure of distance-regular graphs. It can be shown that any sequence of non-increasing $b_i$'s and non-decreasing $c_i$'s give rise to a sequence of potentials $\phi_i$, and that the $\phi_i$'s are decreasing and remain positive. In that sense, a graph doesn't need to actually exist for a given intersection array in order for the potentials to be defined and behave correctly. Furthermore, any intersection arrays which can be ruled out as corresponding to actual graphs by this theorem could in theory be ruled out by the many facts from which we deduced the theorem. Nevertheless, this theorem does perhaps capture a large number of disparate and complicated results on distance-regular graphs in a simple statement. As an example, Theorem \ref{db} shows that the following intersection arrays cannot be realized.

\vski

\begin{tabular}{ l c r }
Intersection array & Vertices & $\frac{\phi_1+ \ldots \phi_{D-1}}{\phi_0}$ \\
\hline
(3,2,2,1,1,1,1;1,1,1,1,1,1,3) & 62 & 1.04918 \\
(5,2,2,1,1,1,1;1,1,1,1,1,1,4) & 101 & 1.0375 \\
(8,3,3,3,3,3,3,3,2,2,1;1,2,2,3,3,3,3,3,3,3,8) & 150 & 0.938852 \\
\end{tabular}

\vski

This can be shown by other methods, but the methods may differ between the examples, and may have much in common with the given proof of Theorem \ref{db} in certain cases. Note that these intersection arrays satisfy a number of basic feasibility requirements, such as being monotone and having $c_i \leq b_{D-i}$ for all $i$. Note further that none of these arrays can be ruled out by Ivanov's bound(Corollary 5.9.6 of \cite{drgraphs}). We therefore have hopes that this theorem can be found useful in the study of distance-regular graphs, both for disallowing certain intersection arrays and as a tool for proving other statements.

\section{Acknowledgements}

The first author was supported by Priority Research Centers Program through the National Research Foundation of Korea (NRF) funded by the Ministry of Education, Science and Technology (Grant \#2009-0094070). The second author was partially supported by the Basic Science Research Program through the National Research
Foundation of Korea(NRF) funded by the Ministry of Education, Science and Technology (Grant \# 2009-0089826).

\end{document}